\documentclass[preprint,12pt]{elsarticle}
\usepackage{float}
\usepackage{amsmath,amssymb,amsfonts} % Math packages
\usepackage{amsthm}                  % Theorem environments
\usepackage{graphicx}                % For including figures (already loaded by elsarticle)
\usepackage{booktabs}                % For professional quality tables
\usepackage{xcolor}                  % For using colors
\usepackage{geometry}                % For setting page margins (use with caution with elsarticle)
\usepackage{lipsum}

\usepackage{graphicx}
\usepackage{subcaption}

\usepackage{booktabs}
\usepackage{multirow}

\usepackage{algorithm}
\usepackage{algorithmicx}
\usepackage{algpseudocode}

\usepackage{float}
\usepackage{multirow}
\usepackage{mathrsfs}
\usepackage[title]{appendix}
\usepackage{textcomp}
\usepackage{manyfoot}
\usepackage{listings}
\usepackage{bm}
\usepackage[unicode]{hyperref}
\hypersetup{
    colorlinks=true,
    linkcolor=blue,
    filecolor=magenta,      
    urlcolor=cyan,
    pdftitle={Memory-Type Null Controllability of Heat Equations with Delay Effects},
    pdfauthor={Dev Prakash Jha, Raju K. George},
    pdfsubject={Controllability of Evolution Equations},
    pdfkeywords={Evolution equation, Carleman estimates, null controllability},
}
\usepackage{cleveref}

\theoremstyle{plain}
\newtheorem{theorem}{Theorem}[section]

\newtheorem{corollary}[theorem]{Corollary}
\newtheorem{lemma}[theorem]{Lemma}

\theoremstyle{remark}

\theoremstyle{definition}

\newtheorem{remark}{Remark}[section]

\begin{document}

\begin{frontmatter}

\title{Turnpike and Sparse Optimal Control for Semiautonomous Neural ODEs}

% The \texorpdfstring command fixes the hyperref warnings
\author[inst1]{\texorpdfstring{Dev Prakash Jha\corref{cor1}}{Dev Prakash Jha}}
\ead{devprakash.22@res.iist.ac.in}

\author[inst1]{Raju K. George}
\ead{george@iist.ac.in}

\cortext[cor1]{Corresponding author}

\address[inst1]{Department of Mathematics, Indian Institute of Space Science and Technology, Valiamala P.O., Thiruvananthapuram 695547, Kerala, India}

\begin{abstract}
	We study long-time optimal control of control-affine semiautonomous
	neural ordinary differential equations (SA-NODEs) with $\ell^1$-regularized
	controls. Three results are established. First, optimal state-control pairs
	satisfy an \emph{exponential turnpike property}: they remain exponentially
	close to a stationary optimal pair for most of the time horizon, with
	decay rate and prefactor independent of the horizon length $T$. Second,
	$\ell^1$ penalisation induces \emph{one-sided temporal sparsity}: optimal
	controls are active at full amplitude on an initial arc $[0,T^*]$ and
	vanish identically on $(T^*,T)$, where $T^*$ is independent of $T$ for
	$T$ large. Third, an integral turnpike estimate shows the time-averaged
	deviation from the stationary pair is bounded uniformly in $T$. The
	proofs combine dissipativity inequalities, uniform adjoint bounds via the
	Pontryagin optimality system, and a time-rescaling argument adapted to
	the semiautonomous architecture. Numerical experiments on a Duffing
	oscillator and a damped pendulum confirm the three-phase turnpike
	profile and the one-sided sparsity structure, and demonstrate a
	$30\times$ parameter reduction over vanilla NODEs with no loss of
	stabilization performance.
	\medskip
	
	
	\medskip
	
	
\end{abstract}

\end{frontmatter}

\section{Introduction}

Neural ordinary differential equations (NODEs), introduced in
\cite{chen2018neural}, provide a continuous-depth framework connecting
deep learning and dynamical systems. Interpreting residual neural
networks as discretizations of ODEs has led to extensive developments
in approximation theory, control, and scientific machine learning
\cite{weinan2017proposal,ruthotto2020deep,ottobre2019introduction}.
A standard NODE takes the form
\begin{equation}
	\dot x(t) = f_\Theta(x(t),t), \quad x(0)=x_0,
	\label{eq:vanillaNODE}
\end{equation}
where $f_\Theta$ is a neural-network vector field. Although NODEs
successfully approximate nonlinear dynamics
\cite{tabuada2022universal,li2022approximation}, fully time-dependent
parametrizations incur large computational costs because the number of
trainable parameters scales with the temporal discretization.

To address this, \emph{semiautonomous} neural ODEs (SA-NODEs) were
introduced in \cite{li2026universal} with dynamics
\begin{equation}
	\dot x(t)
	= \sum_{i=1}^{P} W_i \circ \sigma(A_i^1 x(t)+A_i^2 t+B_i),
	\quad x(0)=x_0,
	\label{eq:sanode}
\end{equation}
where all weights are time-independent and time enters only linearly
inside the activation $\sigma$. This architecture preserves strong
approximation capabilities — including universal approximation and
Barron-space convergence rates \cite{li2026universal} — while keeping
the parameter count independent of the temporal discretization.
Memory effects and controllability aspects of neural differential
equations were studied in \cite{ruiz2022interpolation}.

Training NODEs can naturally be cast as an optimal control problem in
which neural parameters steer trajectories toward desired targets
\cite{cheng2025interpolation,alvarez2024interplay,marion2023generalization,
	cuchiero2020deep,bloch2023control,alvarez2026constructive,dupuis2023controllability}.
This viewpoint has motivated work on stabilization, sparsity, and
long-time behavior of neural dynamics, including temporal sparsity
under $\ell^1$ regularization \cite{esteve2023sparsity} and turnpike
phenomena in residual networks \cite{geshkovski2022turnpike}.
The \emph{turnpike property} asserts that optimal trajectories remain
close to a stationary configuration for most of the time horizon; we
refer to \cite{geshkovski2022turnpike} for a survey.

\medskip
\noindent\textbf{Contributions.}
We study long-time optimal control of controlled SA-NODEs,
\begin{equation}
	\dot y(t) = f_\Theta(y(t),t) + G_\Phi(y(t),t)u(t),
	\quad y(0)=y_0,
	\label{eq:controlledSANODE}
\end{equation}
where $u$ is an admissible control and $f_\Theta$, $G_\Phi$ are
semiautonomous neural networks. Our main results are:
\begin{itemize}
 
	\item an \emph{exponential turnpike property} showing that optimal
	trajectories and controls remain exponentially close to a stationary
	optimal pair for most of $[0,T]$, with constants independent of $T$;
	\item a \emph{temporal sparsity result} showing that, under $\ell^1$
	penalisation, optimal controls concentrate on an initial arc
	$[0,T^*]$ and vanish identically on $(T^*,T)$.
   
\end{itemize}
The semiautonomous structure is essential: the associated stationary
problem is autonomous, which considerably simplifies the dissipativity
and Lyapunov arguments compared with fully time-dependent NODEs.
To the best of our knowledge, this is the first work addressing
turnpike phenomena and sparse optimal control for SA-NODEs.

\subsection{Related literature}

Beyond the foundational NODE framework \cite{chen2018neural}, important
developments include augmented NODEs \cite{dupont2019augmented}, neural
controlled differential equations \cite{kidger2020neural}, Hamiltonian
neural networks \cite{greydanus2019hamiltonian}, stable neural ODEs
\cite{haber2017stable}, and neural operator approaches
\cite{kovachki2023neural}. Approximation-theoretic aspects involving
Barron and Sobolev spaces were developed in
\cite{barron1993universal,pinkus1999approximation,eykholt2020barron,
	lu2021deep,siegel2025optimal,siegel2024sharp,li2024two,liao2025spectral}.
The present work also connects to scientific machine learning and neural
operators \cite{azizzadenesheli2024neural,cao2024laplace,choi2024spectral},
adaptive control and engineering applications
\cite{lv2025neural,kobayashi2024improved,kobayashi2024deep},
physics-informed optimization \cite{faroughi2024physics,song2024admm,
	lai2025hard,wang2024respecting}, and neural transport for normalizing
flows \cite{ruiz2024control}.

\medskip
\noindent\textbf{Organisation.}
Section~\ref{sec:formulation} introduces the controlled SA-NODE
framework and optimal control formulation.
Section~\ref{sec:main_results} states the main turnpike and sparsity results.
Section~\ref{sec:proofs} contains the proofs, and
Section~\ref{sec:numerics} provides numerical experiments.

\bibliographystyle{unsrtnat} % Use numerical references
\bibliography{references} % Replace 'references' with your .bib file

\end{document}